\documentclass[11pt]{amsart}

\usepackage{amsmath}
\usepackage{amssymb}
\usepackage{latexsym}
\usepackage{amscd}
\usepackage{citesort}
\usepackage{mathrsfs}
\usepackage[all]{xy}

\newdimen\AAdi%
\newbox\AAbo%
\def\AAk#1#2{\s_etbox\AAbo=\hbox{#2}\AAdi=\wd\AAbo\kern#1\AAdi{}}%
\def\AAr#1#2#3{\s_etbox\AAbo=\hbox{#2}\AAdi=\ht\AAbo\raise#1\AAdi\hbox{#3}}%
\font\tenmsb=msbm10 at 12pt \font\sevenmsb=msbm7 at 8pt
\font\fivemsb=msbm5 at 6pt
\newfam\msbfam
\textfont\msbfam=\tenmsb \scriptfont\msbfam=\sevenmsb
\scriptscriptfont\msbfam=\fivemsb
\def\Bbb#1{{\tenmsb\fam\msbfam#1}}
\newtheorem{theorem}{Theorem}

\newtheorem{remark}[theorem]{Remark}

\newtheorem{corollary}[theorem]{Corollary}

\newtheorem{lemma}[theorem]{Lemma}
\newtheorem{proposition}[theorem]{Proposition}

\numberwithin{equation}{section} \numberwithin{theorem}{section}

\renewcommand{\topmargin}{0cm}
\renewcommand{\oddsidemargin}{5mm}
\renewcommand{\evensidemargin}{5mm}
\renewcommand{\textwidth}{150mm}
\renewcommand{\textheight}{230mm}

\def\Hess{\mbox{Hess}}

\def\R{\mathbb R}
\def\oR{\overline{\mathbb{R}}}

\def\na{\nabla}
\def\bn{\overline\nabla}

\def\ir#1{\mathbb R^{#1}}

\def\f#1#2{\frac{#1}{#2}}

\def\grs#1#2{\bold G_{#1,#2}}

\def\dd#1#2{\frac {d\,#1}{d\,#2}}
\def\dds#1#2{\frac{d^2\,#1}{d\,#2^2}}
\def\dt#1{\frac {d\,#1}{d\,t}}
\def\mc#1{\mathcal{#1}}

\def\pd#1#2{\frac {\partial #1}{\partial #2}}

\def\U{\Bbb{U}}

\def\a{\alpha}
\def\be{\beta}

\def\r{\Re_{I\!V}}

\def\p#1{\partial #1}

\def\de{\delta}
\def\De{\Delta}
\def\e{\eta}
\def\ep{\varepsilon}
\def\eps{\epsilon}

\def\G{\Gamma}
\def\g{\gamma}
\def\k{\kappa}
\def\la{\lambda}
\def\La{\Lambda}
\def\lan{\langle}
\def\ran{\rangle}
\def\om{\omega}
\def\Om{\Omega}
\def\th{\theta}

\def\w{\wedge}
\def\r{\rho}

\def\div{\mathrm{div}}

\begin{document}

\title[The translating solitons of mean curvature flow]
{The translating solitons of mean curvature flow }

\author{Y.L. Xin}
\address{Institute of Mathematics, Fudan University, Shanghai 200433, China}
\email{ylxin@fudan.edu.cn}
\thanks{The author is supported partially by NSFC}

\begin{abstract}
We study some basic problems of translating solitons: the volume growth, generalized maximum principle, Gauss maps and certain functions
related to the Gauss maps, finally we carry out point-wise estimates and integral estimates for the squared norm of the second fundamental form. Those estimates give rigidity theorems for translating solitons in the Euclidean space in higher codimension.

\end{abstract}

\maketitle

\section{Introduction}

Let $X : M^n \rightarrow \R^{m+n}$ be an isometric immersion from
an $n$-dimensional oriented Riemannian manifold $M$ to the Euclidean space $\R^{n+m}$. Let us consider the mean curvature flow (abbreviated by MCF) for a submanifold  $M$ in
$\ir{m+n}.$ Namely, considering a one-parameter family $X_t=X(\cdot,
t)$ of immersions $X_t:M\to \ir{m+n}$ with corresponding images
$M_t=X_t(M)$ such that
\begin{equation}\left\{\begin{split}
\dt{}X(x, t)&=H(x, t),\qquad x\in M\\
X(x, 0)&=X(x)
\end{split}\right.\label{mcfe}
\end{equation}
is satisfied, where $H(x, t)$ is the mean curvature vector of $M_t$ at $X(x, t)$ in $\ir{m+n}.$

There is a special class of solutions to (\ref{mcfe}), called  \emph{translating solitons} (abbreviated by \emph{translators}) of MCF. A submanifold
$X: M^n\to \ir{m+n}$ is said to be a  translator, if it satisfies
\begin{equation}\label{te}
H= V^{N}.
\end{equation}
Here, $V$ is a constant vector with unit length in $\R^{m+n},$ and $V^N$ is the normal projection of $V$ to the normal bundle of $M$ in $\R^{m+n}$. The translator gives raise an eternal solution $X_t=X+tV$ to MCF equations (\ref{mcfe}).   The translators  play an important role in the study of the mean curvature flow. They are not only special solutions to
the mean curvature flow equations, but they often occur as Type-II singularity of a mean curvature flow \cite{A-V1}, \cite{A-V2}, \cite{H-S}, \cite{Wh1} and \cite{Wh2}. Recently, the subject has been paid considerable attentions, see \cite{CSS},  \cite{Wang}, \cite{Ng1}, \cite{Ng2}, \cite{Hal} and \cite{MSS}.

In $\R^{m+n}$ define a conformally flat
metric $ds^2=e^{\f{2}{m+n}\lan V,X\ran}\left((dx^1)^2+\cdots+(dx^{m+n})^2\right)$  and obtain a Riemannian manifold, which is denoted by $\oR^{m+n}.$
We define a weighted  volume functional $F_{II}$
on the submanifold $M$ in  $\R^{m+n}$ by
\begin{equation}\label{vo}
F_{II}(M)=\int_Me^{\left<V, X\right>}d\mu,
\end{equation}
where $d\mu$ is the volume element of $M$, induced  from the ambient Euclidean metric. $F_{II}(M)$  is also the volume functional of $M$ in $\oR^{m+n}.$
It is easily seen from the first variational formula (see Section 1.2 of \cite{X}, for example) that any translator satisfies the Euler-Lagrangian
equation of the weighted volume functional $F_{II}$. Namely, $M$ is also a minimal submanifold in $\oR^{m+n}.$

We introduce a linear operator on $M$
\begin{equation}\label{L2}
\mc{L}_{II}=\De+\langle
V,\na(\cdot)\rangle=e^{-\left<V,X\right>}\div(e^{\left<V,X\right>}\na(\cdot)),
\end{equation}
in a similar manner of the drift-Laplacian on the self-shrinkers by Colding and Minicozzi in \cite{CM1}. It can be shown that  $\mathcal{L}_{II}$ is self-adjoint respect to the measure $e^{\left<V,X\right>}d\mu.$

In the present paper we study several basic issues for translators in Euclidean space. In the next section we will derive the Bochner type inequalities for the squared norm of the second fundamental form $|B|^2$ for translators in  any codimension in our terminology (Proposition \ref{Sf}).

The translating solitons, as solutions to the variational problem, share the second variational formula.
In the section 3 we introduce second variational formula for oriented translators of codimension one and then define the stability notion (\ref{ja})(\ref{st}), although more general formula holds. Hence, the notion of $L-$stability in \cite{Sh} is just the usual stability notion for oriented minimal hypersurfaces.

Then, we consider the volume growth in the $4^{th}$ section. We can show that any complete translating soliton has infinite volume and has Euclidean volume growth at least (Theorem \ref{vol}). This is in sharp contrast to self-shrinkers \cite{DX}. We also show that graphic translating soliton hypersurfaces are weighted area-minimizing (Theorem \ref{am}), as a corollary they are stable, which generalizes a result in \cite{Sh}.

It is well known that there is no compact translating soliton.  In section 4 by using the L. V. Ahlfors trick we could prove that one can always use generalized maximum principle for any  $C^2-$function with sub-linear growth on a translating soliton (see Theorem \ref{max}), namely, we show that the Omori-Yau  maximum principle is always applicable in the translating soliton situation. We expect its more applications in the further investigation for translators.

For a submanifold in Euclidean space we can define the Gauss map. In the section 6 we firstly prove translating soliton version of the Ruh-Vilms theorem \cite{RV} (see Theorem \ref{RVtype}), as a conclusion that for any translating soliton its Gauss map is a weighted harmonic map (see Corollary \ref{RVtype1}). The image of the Gauss map is a Grassmannian manifold. In our previous work we studied some natural functions ($w-$functions and $v-$functions) in the Grassmannian manifolds (see \cite{XY1} and \cite{JXY1}). Using certain estimates for those functions we can obtain accurate estimates for $\mc{L}_{II}v$ and $\mc{L}_{II}w$ on translating solitons (see Proposition \ref{propw} and Proposition \ref{propv}).

In the final section we derive integral estimates for the squared norm of the second fundamental form $|B|^2$ in terms of the Sobolev constant by using the Sobolev inequality and the Bochner type inequalities derived in the previous section. This leads to a rigidity result (Theorem \ref{sob}). We can also carry out local estimates for $|B|^2$ in terms of the $v-$ function. If the $v-$function satisfies
$$v\le v_1<v_0=\f{2\cdot 3^{\f{2}{3}}}{1+3^{\f{2}{3}}},$$
we obtain  point-wise estimates for $|B|^2$, as shown  by (\ref{geB}). This leads a rigidity result (see first conclusion of Theorem \ref{ge}). If we assume $|B|^2$ is constant by using our maximum principe in the previous section we can relax the conditions on the $v-$function. This is the second conclusion of Theorem \ref{ge}. In general, $v\ge 1$. The upper bound of the $v-$function implies the restriction of the image under the Gauss map. Such a bound is necessary for Bernstein type results of minimal submanifolds in the Euclidean space in higher codimension, as was shown by Lawson-Osserman's counter examples \cite{L-O}.

For a complete translator $M$, if a $v-$function has upper bound, then $M$ could be described by an entire graph $u^\a: \R^n\to\R^m$. The induced metric on $M$ is $ds^2=g_{ij}dx^idx^j$ with $g_{ij}=\de_{ij}+\sum_\a u^\a_iu^\a_j$. Now, the equations
(\ref{te}) becomes
\begin{equation}\label{gte}
g^{ij}u_{ij}^\a=V_\a-\sum_iu^\a_iV_i,
\end{equation}
where $(V_i, V_\a)$ denotes the given constant unit vector. This is an elliptic system. In this case the $v-$function is just
$\sqrt{\det(\de_{ij}+\sum_\a u^\a_iu^\a_j)}.$
Our results claim that  entire solutions to (\ref{gte}) have to be affine linear functions provided
$$\sqrt{\det(\de_{ij}+\sum_\a u^\a_iu^\a_j)}\le v_1<v_0=\f{2\cdot 3^{\f{2}{3}}}{1+3^{\f{2}{3}}}$$.

\bigskip

\section{Bochner type inequalities for the translating soliton}
\medskip

In this section we derive the Bochner type inequality  for
translators in any codimension in our notations. In the case of codimension one the formula is already known (see \cite{MSS} for example).

Let $\na$ and $\overline{\na}$ be the Levi-Civita
connections on $M$ and $\R^{m+n}$, respectively.   Then we define
the second fundamental form $B$ by
$B(V,W)=(\overline{\na}_VW)^N=\overline{\na}_VW-\na_VW$ for any tangent vector fields
$V,W$ in $M$, where $(\cdots)^N$ stands for the orthogonal
projection into the normal bundle $NM$. Similarly, $(\cdots)^T$ stands for
the tangential projection.   The mean curvature vector
$H$ of $M$ is given by
$H=\mathrm{trace}(B)=\sum_{i=1}^nB(e_i,e_i)\in \G(NM),$ where $\{e_i\}$ is
a local orthonormal frame field of $M$. In what follows we use $\na$ for natural
connections on various bundles for notational simplicity if there
is no ambiguity from the context. For $\nu\in\G(NM)$ the shape
operator $A^\nu: TM\to TM$, defined by $A^\nu(V)=-(\bn_V\nu)^T$,
satisfies $\left<B_{V W}, \nu\right>=\left<A^\nu(V), W\right>.$

For minimal submanifolds in an arbitrary ambient Riemannian
manifold J. Simons \cite{Si} derived the Laplacian of the squared
norm of the second fundamental form. For arbitrary submanifolds
in Euclidean space Simons type formula was also  derived (see
\cite{Si}, \cite{X1}, for example).

Choose a local orthonormal frame field $\{e_i,e_\a\}$ along $M$
with dual frame field $\{\omega_i,\omega_\a\}$, such that $e_i$
are tangent vectors of $M$ and $e_{\a}$ are normal to $M$. The
induced Riemannian metric of $M$ is given by $ds_M^2
=\sum_i\omega_i^2$ and the induced structure equations of $M$ are
\begin{equation*}\begin{split}
   & d\omega_i = \omega_{ij}\wedge\omega_j,\qquad
                 \omega_{ij}+\omega_{ji} = 0,\cr
   &d\omega_{ij}= \omega_{ik}\wedge\omega_{kj}+\omega_{i\a}\wedge\omega_{\a j},\cr
   &\Omega_{ij} = d\omega_{ij}-\omega_{ik}\wedge\omega_{kj}
                = -\frac 12 R_{ijkl}\omega_k\wedge\omega_l.
\end{split}\end{equation*}
By Cartan's lemma we have
$$\omega_{\a i} = h_{\a, ij}\omega_j.$$
Here and in the sequel we agree with the following range of
indices
$$1\le i, j, k, \cdots \le n,\quad n+1\le \a, \be, \g, \cdots \le n+m.$$

Set
$$B_{ij}=B_{e_ie_j}=(\bar\na_{e_i}e_j)^N=h_{\a, ij}e_\a,\quad S_{\a\be}=h_{\a, ij}h_{\be, ij}.$$
Then, $$|B|^2=\sum_{\a}S_{\a\a}.$$

From Proposition 2.2 in \cite{X1} we have
\begin{equation}\label{LB}\begin{split}
\De|B|^2=2\,|\na B|^2+2\,\left<\na_i\na_jH,B_{ij}\right>
       &+\,2\left<B_{ij},H\right>\left<B_{ik},B_{jk}\right>\\
       &\quad -\,2\,\sum_{\a\ne\be}|[A^{e_\a}, A^{e_\be}]|^2-2\,\sum_{\a,\be}S_{\a\be}^2.
\end{split}\end{equation}

Suppose that the above chosen frame field
$\{e_i\,\, e_\a\}$ is normal at a considered point $p\in M$. From
the translator equations (\ref{te}) we obtain
\begin{equation}\label{DH}
\na_jH=-\left<V, e_k\right>B_{jk},
\end{equation}
and
\begin{equation}\label{D2H}
\na_i\na_j H=-\langle H, B_{ik}\rangle
B_{jk}-\langle V,e_k\rangle\na_i B_{jk}.
\end{equation}
Combining (\ref{te}), (\ref{L2}), (\ref{LB}) and (\ref{D2H}) (and using the Codazzi  equation), we have
\begin{equation}\aligned\label{L2B}
\mc{L}_{II}|B|^2=2|\na B|^2-\,2\,\sum_{\a\ne\be}|[A^{e_\a},
A^{e_\be}]|^2-2\,\sum_{\a,\be}S_{\a\be}^2.
\endaligned
\end{equation}
This is the translator version of the well-known Simons'
identity. In particular, when the codimension $m=1$ or the normal bundle is flat, the above
Simons' type identity reduces to the following one:

\begin{equation*}\label{LB1}
\mc{L}_{II}|B|^2\ge 2|\na B|^2 - 2|B|^4.
\end{equation*}

When the codimension $m\ge 2$ there is estimate
\cite{LL}\cite{ChenX}
$$\sum_{\a\ne\be}|[A^{e_\a}, A^{e_\be}]|^2+\sum_{\a, \be}S_{\a\be}^2\le\frac 32 |B|^4.$$
Combining (\ref{L2B}) and the above inequality, we have
\begin{equation*}\aligned\label{LB2}
\mathcal{L}_{II}|B|^2\ge 2|\na B|^2-3|B|^4.
\endaligned
\end{equation*}

We summarize the above calculations as the following proposition.

\begin{proposition}\label{Sf}
The squared norm of the second fundamental form $|B|^2$ of the translators in Euclidean space $\R^{m+n}$ satisfies  Simons' type inequality:
\begin{equation}\label{LB1}
\mc{L}_{II}|B|^2\ge 2|\na |B||^2 - k\,|B|^4,
\end{equation}
where
\begin{eqnarray}\label{k}
   k= \left\{\begin{array}{ccc}
     2      \quad\ \ \ {\rm{if}}\ \ \ & m=1 & {\rm{ or\ the\ normal\ bundle\ is\ flat,}} \ \ \  \\
     3      \quad\ \ \ {\rm{if}} \ \ \ &  m\ge 2, &
     \end{array}\right.
\end{eqnarray}
\end{proposition}

\bigskip

\section{Second variational formula for translation solitons}

\medskip
The translators are also the solutions to a natural geometric variational problem. The equation (\ref{te}) is the Euler-Lagrangian Equation of the weighted
volume functional (\ref{vo}). We can also derive the second variational formula from the minimal surface theory. Here, we only consider the oriented codimension one case for possible applications.

From the first variational formula of the volume functional we have
$$\dd{}{s}F_{II}=\int_M\left<V-H_s, \dd{f_s}{s}\right>e^{\left<V,X_s\right>}d\mu_s$$
for a smooth one-parametric family of immersions $f_s$ of $M$ into $\R^{m+n}$ with $f_0=M$. When $M$ is an oriented hypersurface choose
the variational vector field
$$\left.\dd{f_s}{s}\right|_{s=0}=\phi\nu,$$
where $\nu$ is the unit normal vector field of $M$ in $\R^{m+n}$ and $\phi$ is any smooth function on $M$ with
compact support. Hence, the above first variational formula becomes
\begin{equation}\label{fv}
\dd{}{s}F_{II}=\int_M\left<V-H_s, \phi\nu_s\right>e^{\left<V,X_s\right>}d\mu_s
\end{equation}
with unit normal vector field $\nu_s$ to $f_s(M)$ in $\R^{n+1}$. From
$$\left<(X_s)_i, \nu_s\right>=0$$
we obtain
$$\left<\left(\dd{X_s}{s}\right)_i,\nu_s\right>+\left<(X_s)_i, \dd{\nu_s}{s}\right>=0.$$
(Here, $(\cdot)_i$ denote $\na_{e_i}(\cdot)$ with respect to a local orhtonormal frame field $\{e_j\}$ in $M$.)
At $s=0$
$$\left<\phi_i\nu+\phi\nu_i,\nu\right>+\left<X_i,\left.\dd{\nu_s}{s}\right|_{s=0}\right>=0$$
and then
$$\left<X_i,\left.\dd{\nu_s}{s}\right|_{s=0}\right>=-\phi_i.$$
It follows that
\begin{equation}\label{dnut}
\left.\dd{\nu_s}{s}\right|_{s=0}=-\na\phi.
\end{equation}
On the other hand, from the derivation of the second variational formula of minimal submanifolds (see pp 146, 147 and 154 in \cite{X}) we have
\begin{equation}\label{dH}
\left.\dd{}{s}H_s\right|_{s=0}=(\De\phi+|B|^2\phi)\nu.
\end{equation}

Now, let us derive the second variational formula for oriented translating soliton hypersurfaces in the Euclidean space. From (\ref{te}) and (\ref{fv}) we have
\begin{equation}
\dds{}{s}F_{II}(0)=\int_M\left.\dd{}{s}\left<V-H_s, \phi\nu_s\right>\right|_{s=0}e^{\left<V,X\right>}d\mu.
\end{equation}
Noting (\ref{dnut}) and (\ref{dH}), the above expression becomes
\begin{equation}\label{sv}
\dds{}{s}F_{II}(0)=-\int_M\phi(\De\phi+|B|^2\phi+\left<V,\na\phi\right>d\mu=-\int_M\phi(\mc{L}_{II}\phi+|B|^2\phi)d\mu.
\end{equation}
This is the second variational formula for oriented translating solitons of codimension one in Euclidean space. The associate Jacobi operator is defined by
\begin{equation}\label{ja}
J=-\mc{L}_{II}-|B|^2.
\end{equation}

For an oriented  translating soliton $M$ in $\R^{n+1}$ of codimension one, as a minimal hypersurface in the conformally flat $\oR^{n+1}$, we could also have stability notion. If for any smooth function $\phi$ with compact support
\begin{equation}\label{st}
\int_M\phi\,J\phi\,d\mu\ge 0
\end{equation}
is satisfied, $M$ is called a stable translation soliton.

\bigskip

\section{Volume growth for translators}

\medskip

Take orthonormal basis
$$X_\a=e^{-\f{1}{n}x^n}\pd{}{x^\a}, \quad i=1,\cdots n-1,\quad \a=1, \cdots, n$$
in $\oR^{n}$, $\R^n$ with conformally flat metric $ds^2=e^{\f{2}{n}x^n}\left((dx^1)^2+\cdots+(dx^n)^2\right).$
By  direct computations we have
$$[X_i, X_j]=0,\qquad [X_i,X_n]=\f{1}{n}e^{-\f{2}{n}x^n}\pd{}{x^i}$$
and
$$\na_{X_i}X_i=-\f{1}{n}e^{-\f{2}{n}x^n}\pd{}{x^n}, \quad  \na_{X_i}X_j=0, \quad  \na_{X_n}X_i=0, \quad \na_{X_n}X_n=\f{1}{n}\pd{}{x^n}.$$
It follows that
$$R_{ijij}=\lan -\na_{X_i}\na_{X_j}X_i+\na_{X_j}\na_{X_i}X_i+\na_{[X_i,X_j]}X_i, X_j\ran=-\f{1}{n^2}e^{-\f{2}{n}x_n},$$  $$R_{ini}=-\na_{X_i}\na_{X_n}X_i+\na_{X_n}\na_{X_i}X_i+\na_{[X_i,X_n]}X_i=0.$$

We summarize the following conclusion.

\begin{lemma}
Let $V$ be a constant vector of unit length in the Euclidean space $\R^n$.  When it is endowed with a conformally flat metric with the conformal factor
$e^{\f{2}{n}\lan V, X\ran}$, where $X$ stands for the position vector in $\R^n$. Then the sectional curvature
$$\k=\left\{\begin{array}{ccc}
     0      \ \ \qquad\qquad  {\rm{if\ the\ section\ containing}}\ \ V, \  \ \ \  \\
     -\f{1}{n^2}e^{-\f{2}{n}\lan V,X\ran}      \quad {\rm{if\ the\ section\ perpendicular\ to}} \ V.\ \  &
     \end{array}\right.$$
Hence, in $\oR^{n}$ the sectional curvature satisfies $-\f{1}{n^2}e^{-\f{2}{n}\lan V,X\ran}\le \k\le 0.$
\end{lemma}

Since $\De\lan V, X\ran=|H|^2\ge 0$ for any translator there is no maximum of $\lan V, X\ran$ and an similar argument as that in the section 4 of \cite{X2} shows that the volume growth of a translator is Euclidean volume growth at least, as the following consideration.

Let $M$ be a translator in $\R^{m+n}$, as well as a minimal submanifold in $\oR^{m+n}$. Fix a point $o\in M\subset \oR^{m+n},$ denote the distance function from  $o$ in $\oR^{m+n}$ by $\r.$ By using the classical Hessian comparison theorem for  $\oR^{m+n}$ with non-positive sectional curvature, we have
\begin{equation}\label{her}
\overline{\text{Hess}}(\r)(X, Y)\ge \f{1}{\r}(\left<X, Y\right>-\left<X,\na\r\right>\left<Y,\na\r\right>)
\end{equation}
for $X, Y\in T\oR^{m+n}.$ The restriction of $\r$ on $M$ gives that for $X, Y\in TM\subset T\oR^{m+n}.$
\begin{equation}\label{her1}\aligned
\text{Hess}(\r)(X, Y)&=XY(\r)-(\bn_{X}Y-(\bn_XY)^N)\r\\
&=\overline{\text{Hess}}(\r)(X, Y)+\left<\overline{B}_{XY},\na\r\right>,
\endaligned
\end{equation}
where $\overline{B}$ stands for the second fundamental form of $M$ in $\oR^{m+n}$
It follows that
\begin{equation}\label{Lr}\aligned
\De\r&\ge\f{1}{\r}(n-|\na\r|^2_M)+\left<\bar H,\na\r\right>,\\
\De\r^2&=2|\na\r|^2_M+2\r\De\r\ge 2n+2\r\left<\bar H,\na\r\right>\ge 2n,
\endaligned\end{equation}
since $M$ is minimal in $\oR^{m+n}$ and its mean curvature $\bar H=0$.  Let $B(\r)$ be a geodesic ball of radius $\r$ and centered at $o\in M\subset\oR^{m+n}$. Its restriction on $M$ is denoted
by
$$D(\r)=B(\r)\bigcap M.$$
Obviously, $\pd{}{\r}$ is the unit normal vector to $\partial B(\r)$. Its orthogonal projection to $M$ is normal to $\partial D(\r)$. Let $\{e_\a\}$ be
orthonomal normal frame field of $M$ in $\oR^{m+n}$ at the considered point. Then
$$\pd{}{\r}-\left<\pd{}{\r},e_\a\right>e_\a=\pd{}{\r}-e_a(\r)e_\a$$
is normal to $\partial D(\r)$. Hence,
$$\nu=a\left(\pd{}{\r}-e_a(\r)e_\a\right)$$
is the unit normal vector to $\partial D(\r)$, where $a=\f{1}{\sqrt{1-\sum_\a(e_\a(\r))^2}}$. Noting
$$\aligned
\left<\pd{}{\r},\na\r^2\right>&=2\r,\\
\left<e_\a,\na\r^2\right>&=2 \r\, e_\a(\r),
\endaligned$$
and then
$$\aligned\left<\nu,\na\r^2\right>&=a(2\r-\sum_\a e_\a(\r)\cdot 2\r\, e_\a(\r))\\
&=2\r\sqrt{1-\sum_\a(e_\a(\r))^2}\le 2\r.\endaligned$$
Integrating (\ref{Lr}) over $D(\r)$ and using Stokes' theorem, we have
\begin{equation}\label{ILr}\aligned
2\, n\,\text{vol}(D(\r))&\le\int_{D(\r)}\De\r^2*1\\
&=\int_{D(\r)}\left<\nu,\na\r^2\right> *1\\
&\le 2\,\int_{\partial D(\r)}\r *1=2\,\r\, \text{vol}(\partial D(\r)).
\endaligned\end{equation}
On the other hand, since $|\na\r|_M\le 1,$ the co-area formula gives
\begin{equation}\label{co}
\f{d\text{vol}(D(\r)))}{d\r}\ge\text{vol}(\partial D(\r)).
\end{equation}
Combining (\ref{ILr}) and (\ref{co}), we obtain that
$$\f{d(\text{vol}(D(\r)))}{\text{vol}D(\r)}\ge\f{n\,d\r}{\r}.$$
It follows that
\begin{equation}\label{v}
\f{\text{vol}(D(\r))}{\r^n}
\end{equation}
is a nondecreasing function in $\r$, which implies vol$(M)$ is infinite and it has polynomial growth in $\r$ of order $n$ at least. We call such minimum volume growth as the Euclidean volume growth (note that $\r$ is distance function in $\oR^{m+n}$) . We conclude the following results.

\begin{theorem}\label{vol}
Any complete translating soliton in the Euclidean space has infinite volume and has Euclidean volume growth at least.
\end{theorem}

Let us consider the graphic translator of codimension one situation. Choose coordinates so that the given vector $V$ along the $x^{n+1}$ axis and the hypersurface $M$ is defined by $(x, u(x))$. Now, the translator equation (\ref{te}) becomes
\begin{equation}\label{te1}
\div \left(\f{Du}{\sqrt{1+|Du|^2}}\right)=\f{1}{\sqrt{1+|Du|^2}}.
\end{equation}

In \cite{Sh}, the author studied graphic translator surfaces in $\R^3$ and showed it is stable. This is crucial for his subsequent argument. In fact, we can
prove a more general result. A graphic translator hypersurface, as a minimal graphic hypersurface in $\oR^{n+1}$, is area-minimizing by an easy arguments, as similar to minimal graphic hypersurfaces   in $\R^{n+1}$ (see section 6.2 in \cite{X}, for example).

\begin{theorem}\label{am}
Let $\Om$ be a bounded domain in $\R^n$ and $M$ be a  graphic translator on $\overline{\Om}$  with volume element $d\mu_M$. For any hypersurface $W$ in $\R^{n+1}$ with $\p M=\p W$, one has
\begin{equation}\aligned
\int_{M}e^{x^{n+1}}d\mu_M\le\int_{W}e^{x^{n+1}}d\mu_{W},
\endaligned
\end{equation}
where the above inequality attains equality if and only if $W=M$.
\end{theorem}
\begin{proof}
Put $E_i=\pd{}{x^i}, \; E_{n+1}=\pd{}{x^{n=1}}$ and $v=\sqrt{1+|Du|^2}.$ Note that $u_i$ and $v$ are functions on $\Om$. Let $U$ be the domain in $\R^{n+1}$ enclosed by $M$ and $W$. Let $Y$ be a vector field in $U\subset\R^{n+1}$ defined by
$$Y=\left(-\f{u_i}{v}E_i+\f{1}{v}E_n\right)e^{x^{n+1}}.$$
From translator equation (\ref{te1}) we have
\begin{equation*}
\overline{\div}(Y)=-e^{x^{n+1}}\sum_i\pd{}{x^i}\left(\f{u_i}{v}\right)+\f{1}{v}e^{x^{n+1}}=0,
\end{equation*}
where $\overline{\div}$ stands for the divergence operator on $\R^n$.
Let $\nu_M,\nu_W$ be the unit outside normal vectors of $M,W$ respectively. Observe that $Y|_M=\nu_Me^{x^{n+1}}$. Then by Green's formula, up to a minus sign, we have
\begin{equation*}\aligned\label{Gaussf}
0=&\int_{U}\overline{\div}(Y)=\int_{M}\lan Y,\nu_M\ran d\mu_M-\int_{W}\lan Y,\nu_W\ran d\mu_{W}\\
\ge&\int_{M} e^{x^{n+1}}d\mu_M-\int_{W}e^{x^{n+1}} d\mu_{W}.
\endaligned
\end{equation*}
Obviously, the above inequality attains equality if and only if $M=W$.

The above argument is for the situation that hypersurface $M$ lies above the hypersurface $W$. For the opposite case the conclusion is also valid by a similar
argument. Combine these two situations, we conclude that the final inequality always valid for any complicated configuration.
\end{proof}

\bigskip

\section{Maximum principle for translators}

\medskip
Since the height function in the given $V$ direction has no local maximum, there is no compact translating soliton. It is natural to consider the generalized
maximum principle in translators. The following results might be a useful analytic tool for further investigation for translating solitons. In the final section we will give an application to obtain a rigidity result (see second part of the Theorem \ref{ge}).

\begin{theorem}\label{max}
Let $X:M\to \R^{m+n}$ be an $n-$dimensional complete translator. Let $r(x)$
be extrinsic distance function from a fixed point in $M$, namely the restriction to $M$ of the ambient Euclidean distance $|X|$. Then for any $C^2$-function $f$ on $M$, satisfying
$$\lim_{r\rightarrow+\infty}\frac{f(x)}{r}=0,$$
there exists $\{x_k\},$ such that
$$\lim_{k\rightarrow\infty}f(x_k)=\sup f,\ \ \ \lim_{k\rightarrow\infty}|\nabla f(x_k)|=0
\ \ \ \lim_{k\rightarrow\infty}\mc{L}_{II} f(x_k)\leq0.$$
\end{theorem}

\begin{proof}
First of all, from the translator equation (\ref{te})
\begin{equation}\label{dr}\aligned
\na r^2&=2X^T,\quad\text{and}\quad |\na r|^2\le 1,\\
\De r^2&=2n+2\lan H, X\ran=2n+2\lan V^N,X\ran\le 2n+2r,\\
\De r&\le \f{n}{r}+1.
\endaligned
\end{equation}

Take  a sequence of positive numbers $\{\epsilon_k\}$ and  $\epsilon_k\rightarrow 0$ when $k\rightarrow\infty.$ Put
\begin{equation}\label{fk}
f_k(x)=f(x)-\epsilon_kr(x),
\end{equation}
where $f_k(x)$ is  a $C^1-$ function on $M$ and is also a $C^2-$function except the origin.
Since the assumption on $f$, we know that for each $k$
$$f_k\rightarrow-\infty,\qquad \text{as} \quad r\rightarrow+\infty.$$
It follows that there exists a point $x_k$, which is the maximum
point of $f_k$:
$$f_k(x_k)=\sup f_k(x).$$
We then have for any point $x\in M$
$$
f_k(x_k)=f(x_k)-\epsilon_k r(x_k) \geq f_k(x)=f(x)-\epsilon_k r(x),
$$
and
$$f(x_k)\geq f(x)+\epsilon_k\Big(r(x_k)-r(x)\Big).$$
If $r(x_k)$ is bounded
$$\lim_{k\rightarrow\infty}\epsilon_k(r(x_k))-r(x))=0,$$
otherwise, for sufficiently large $k$,
$f(x_k)>f(x)$. Those facts implies that
$$\sup_k f(x_k)\ge\sup_{x\in M} f(x),$$
and in turn
$$\sup_k f(x_k)=\sup_M f(x).$$
Up to a subsequence we have
$$\lim_{k\to\infty}f(x_k)=\sup f(x).$$
This proves the first claim of the Theorem.

If $r(x_k)$ is bounded, there is a convergent subsequence $x_{i_k}$, still denoted by $x_k$ and $f(x)$ attains its maximum at a point in $M$. Otherwise,
$r(x_k)$ are bounded away from the origin and $f_k(x)$ is $C^2-$function near each point $x_k$. Hence,
\begin{equation}\label{mx}
\nabla f_k(x_k)=0,\ \ \ \ \ \ \ \Delta f_k(x_k)\leq 0.
\end{equation}
From (\ref{dr}), (\ref{fk}) and (\ref{mx}) we have
\begin{equation*}
\lim_{k\rightarrow\infty}|\nabla f(x_k)|
=\lim_{k\rightarrow\infty}\epsilon_k|\nabla r(x_k)|=0,
\end{equation*}
\begin{equation*}\lim_{k\rightarrow\infty}\Delta f(x_k)
 \leq\lim_{k\rightarrow\infty}\epsilon_k\Delta r(x_k)
 \leq\lim_{k\rightarrow\infty}\epsilon_k\left(\f{n}{r(x_k)}+1\right)=0.
\end{equation*}
From (\ref{L2}) and the above inequalities give us
$$\lim_{k\rightarrow\infty}\mc{L}_{II} f(x_k)\le 0.$$
\end{proof}

\bigskip

\section{Gauss maps for translators}

\medskip

If $M$ is an oriented submanifold in $\R^{m+n}$, we can define the
Gauss map $\g: M\to \grs{n}{m}$ that is obtained by parallel
translation of $T_pM$ to the origin in the ambient space $\ir{m+n}.$
Here $\grs{n}{m}$ is the Grassmannian manifolds constituting of all oriented $n$-subspaces in
 $\R^{m+n}$. In particular when the codimension $m=1$, $\grs{n}{1}$ becomes
Euclidean sphere.

Using Pl\"ucker coordinates, the Gauss map $\g$ could be described as
$\g(p)=e_1\w\cdots\w e_n$. There are relations (see  (3.2) and (3.4) in \cite{JXY1})

\begin{equation}\aligned\label{ghij}
\g^*\om_{i\a}&=h_{\a, ij}\om_j,\\
|d\g|^2&=\sum_{\a,i,j}h_{\a, ij}^2=|B|^2.
\endaligned
\end{equation}

Now, we consider the Gauss map for  translators in the Euclidean space. There is a notion, so-called $f-$harmonic maps. It was introduced in \cite{Li}. See also \cite{CJQ} for further generalizations. Let $u:M\to N$ be a smooth map. For a smooth function $f$ on $M$. If $u$ is a critical point of the functional
\begin{equation*}
\int_M |du|^2e^f\, d\mu,
\end{equation*}
then $u$ is called an $f-$harmonic map. The Euler-Lagrangian equation of the above $f-$functional is
\begin{equation}\label{ften}
\tau_f(u)=\tau(u)+du(\na f),
\end{equation}
where $\tau(u)=(\na_{e_i} d u)e_i$ is the tension field of the given map $u$, where $\{e_i\}$ is a local orthonormal frame field in $M$. We use the same notation as in \cite{X'}.

By simple calculations we show that the Gauss map for a translator is an $f-$harmonic map with $f=\left<V,X\right>$.

\begin{theorem}\label{RVtype}
For an oriented $n-$dimensional submanifold $X:M\to\R^{m+n}$ its
Gauss map $\g: M\to \grs{n}{m}$ is $f-$ harmonic map with $f=\left<V,X\right>$ if
and only if $H-V^N$ is a parallel vector field in the normal
bundle $NM$, where $V$ is a constant vector of unit length in $\R^{m+n}$.
\end{theorem}

\begin{proof}

Let $\{e_1,\cdots,e_n\}$ be a local tangent orthonormal frame field
on $M$ and $\{e_{n+1},\cdots, e_{m+n}\}$ be a local normal orthornormal
frame field on $M$, and we assume $\na e_i=0$ and $\na e_\a=0$ at
the considered point. From (\ref{ghij})
\begin{equation}\label{tension1}
d\g(e_i)=h_{\a,ij} e_{\a j}
                \end{equation}
where $\{h_{\a,ij}=\lan B_{e_ie_j},e_\a \ran\}$ are coefficients of the second fundamental form, and
$e_{\a j}$ is obtained by replacing $e_j$ by $e_\a$ in
$e_1\w\cdots\w e_n$. We note that $\{e_{\a j}\}$ is an orthornormal basis of the tangent space of
$\grs{n}{m}$ at $e_1\w \cdots\w e_n$.

At the considered point,
\begin{equation}\aligned
\na_{e_i}e_{\a j}=&\na_{e_i}(e_1\w \cdots\w e_\a\w \cdots\w e_n)\\
                =&\sum_{k}e_1\w\cdots\w \na_{e_i}e_k\w\cdots \w e_\a\w\cdots\w e_n\\
                &\qquad+e_1\w \cdots\w \na_{e_i}e_\a\w \cdots\w e_n=0.
                \endaligned
\end{equation}

Using the Codazzi equations one can obtain
\begin{equation}\aligned\label{tension2}
\na_{e_i}h_{\a,ij}&=\na_{e_i}\lan B_{e_i e_j},e_\a\ran=\lan (\na_{e_i}B)_{e_ie_j},e_\a\ran\\
                 &=\lan (\na_{e_j}B)_{e_ie_i},e_\a\ran= \na_{e_j}H^\a
                 \endaligned
\end{equation}
with $H^\a:=\lan H,e_\a\ran$ the coefficients of the mean
curvature vector.

Combining with (\ref{tension1})-(\ref{tension2}) gives
\begin{equation}\label{tension3}\aligned
(\na_{e_i}d\g)e_i&=\na_{e_i}d\g(e_i)=(\na_{e_i}h_{\a,ij})e_{\a j}+h_{\a,ij} \na_{e_i}e_{\a j}\\
&=(\na_{e_j}H^\a) e_{\a j}.
\endaligned
\end{equation}

Since in our translator case $f=\left<V, X\right>$,
\begin{equation}\label{tension4}
\na_{e_i}f=\left<V, e_i\right>.
\end{equation}
Let $V^\a:=\lan V,e_\a\ran$, then
\begin{equation}\label{tension5}\aligned
\na_{e_j}V^\a&=\lan V,\bn_{e_j}e_\a\ran=-\lan V, h_{\a, ij}e_i\ran\\
            &=-h_{\a, ij}\lan V,e_i\ran.
            \endaligned
\end{equation}

In conjunction with (\ref{tension1}), (\ref{tension3}), (\ref{tension4}) and (\ref{tension5}) we have
\begin{equation}\aligned
\tau_f(\g):=&\tau(\g)+d\g(\na f)\\
              =&\Big[\na_{e_j}H^\a+\left<V, e_i\right>h_{\a, ij}\Big]e_{\a j}\\
              =& \na_{e_j}(H^\a-V^\a)e_{\a j}.
              \endaligned
              \end{equation}
Now, the conclusion follows.
\end{proof}

\begin{corollary}\label{RVtype1}
If $M$ is a translator in $\R^{m+n},$ then its Gauss map $\g:M\to
\grs{n}{m}$ is a $f-$ harmonic map with $f=\lan V, X\ran$.
\end{corollary}

\begin{remark}
The above results are translator version of the Ruh-Vilms theorem in minimal surface theory. This was already be obtained in codimension $m=1$ case by
\cite{BS}.
\end{remark}

Now we assume $F$ to be a $C^2$-function on $\grs{n}{m}$, then $\phi=F
\circ \g$ gives a $C^2$-function on $M$. We also choose a local
orthonormal frame field $\{e_i, e_\a\}$ along $M$ such that $\na e_i=0$ and
$\na e_\a=0$ at the considered point. Denote $f=\lan V,X\ran$ here and in the sequel. A straightforward calculation
shows
\begin{equation*}\aligned
\mc{L}_{II}\phi&=e^{-f}\div(e^f\na \phi)=e^{-f}\na_{e_i}\big(e^f d\phi(e_i)\big)\\
       &=e^{-f}\na_{e_i}\big(e^f dF(d\g(e_i))\big)=e^{-f}\na_{e_i}\big(dF(e^fd\g(e_i))\big)\\
       &=e^{-f}(\na_{e_i}(dF))(e^f d\g(e_i))+e^{-f}dF\Big(\big(\na_{e_i}(e^f d\g)\big)e_i\Big)\\
       &=\text{Hess}\ F(\g_* e_i,\g_* e_i)+dF(d\g(\na f)+\tau(\g))\\
       &=\text{Hess}\ F(\g_* e_i,\g_* e_i)+dF(\tau_f(\g)).
       \endaligned
       \end{equation*}
If $M$ is a translator, by Corollary \ref{RVtype1}, we have the following composition formula for the translators
\begin{equation}\label{com}
\mc{L}_{II}\phi=\text{Hess}\ F(\g_* e_i,\g_* e_i).
\end{equation}

$\grs{n}{m}$ can be viewed as a submanifold of some  Euclidean space
via the Pl\"ucker embedding. The restriction of the
Euclidean inner product on $M$ is denoted by $w:\grs{n}{m}\times \grs{n}{m}\to \R$
$$w(P,Q)=\lan e_1\w\cdots\w e_n,f_1\w\cdots\w f_n\ran=\det W$$
where $P$ is spanned by a unit $n$-vector $e_1\w\cdots\w e_n$, $Q$ is spanned by another unit $n$-vector $f_1\w\cdots
\w f_n$, and $W=\big(\lan e_i,f_j\ran\big)$. It is well-known that
$$W^T W=O^T \La O$$
with $O$ an orthogonal matrix and
$$\La=\left(\begin{array}{ccc}
            \mu_1^2 &   &  \\
                    & \ddots &  \\
                    &        & \mu_n^2
            \end{array}\right).$$
Here each $0\leq \mu_i^2\leq 1$. Putting $p:=\min\{m,n\}$, then
at most $p$ elements in $\{\mu_1^2,\cdots, \mu_n^2\}$ are not
equal to $1$. Without loss of generality,  we can assume
$\mu_i^2=1$ whenever $i>p$. We also note that the $\mu_i^2$ can
be expressed as
\begin{equation}\label{di1a}
\mu_i^2=\frac{1}{1+\la_i^2}. \end{equation}

The Jordan angles between $P$ and $Q$ are defined by
$$\th_i=\arccos(\mu_i)\qquad 1\leq i\leq p.$$
The distance between $P$ and $Q$ is defined by
\begin{equation}\label{di}
d(P, Q)=\sqrt{\sum\th_i^2}. \end{equation}
Thus, (\ref{di1a}) becomes
\begin{equation}\label{di2}
\la_i=\tan\th_i. \end{equation}

In the sequel, we shall assume $n\geq m$ without loss of generality. Let $\a=n+\a'$ and denote $\a$ for $\a'$ for simplicity.

Now we fix $P_0\in \grs{n}{m}.$ We represent it by the $n-$vector
$\eps_1\wedge\cdots\wedge\eps_i\wedge\cdots\wedge\eps_n$. We choose  $ m $ vectors $ \eps_{n+\a} $,
such that $ \{\eps_i, \eps_{n+\a} \} $ form an orthonormal base of $
\ir{m+n} $. Denote
$$\Bbb{U}\triangleq \{P\in \grs{n}{m},\; w(P,P_0)>0\}.$$
Our fundamental quantity will be
\begin{equation}
v(\cdot, P_0):=w^{-1}(\cdot, P_0)\; \text{ on }\; \Bbb{U}.
\end{equation}
For convenience, we define subsets in $\Bbb{U}$ by
$$\Bbb{U}_2\triangleq \{P\in \Bbb{U}, v(\cdot,P_0)<2\},\qquad  \Bbb{U}_3\triangleq \{P\in \Bbb{U}, v(\cdot,P_0)<3\}.$$
For arbitrary $P\in \U$ determined by an $n\times m$ matrix $Z$, it
is easily seen that
\begin{equation}\label{v}
v(P,P_0)=\big[\det(I_n+ZZ^T)\big]^{\f{1}{2}}=\prod_{\a=1}^m
\sec\th_\a =\prod_{\a=1}^m \frac{1}{\mu_\alpha}.
\end{equation}
where $\th_1,\cdots,\th_m$ denote the Jordan angles between $P$
and $P_0$.

In this terminology, Hess$(v(\cdot, P_0)$ has been estimated in
\cite{XY1}. By (3.8) in \cite{XY1}, we have
\begin{eqnarray}\label{He}\aligned
\Hess(v(\cdot,P_0))&=\sum_{i\neq \a}v\ \om_{i\a}^2+\sum_\a
(1+2\la_\a^2)v\ \om_{\a\a}^2
+\sum_{\a\neq\be} \la_\a\la_\be v(\om_{\a\a}\otimes \om_{\be\be}+\om_{\a\be}\otimes\om_{\be\a})\\
&=\sum_{m+1\leq i\leq n,\a}v\
\om_{i\a}^2+\sum_{\a}(1+2\la_\a^2)v\ \om_{\a\a}^2
                                          +\sum_{\a\neq \be}\la_\a\la_\be v\ \om_{\a\a}\otimes\om_{\be\be}\\
&\qquad\qquad+\sum_{\a<\be}\Big[(1+\la_\a\la_\be)v\Big(\f{\sqrt{2}}{2}(\om_{\a\be}
+\om_{\be\a})\Big)^2\\
&\hskip2in+(1-\la_\a\la_\be)v\Big(\f{\sqrt{2}}{2}(\om_{\a\be}-\om_{\be\a})\Big)^2\Big].
\endaligned
\end{eqnarray}
It follows that
\begin{equation}\label{hess}
 v(\cdot,P_0)^{-1}\Hess(v(\cdot,P_0))
 =g+\sum_\a 2\la_\a^2 \om_{\a\a}^2+\sum_{\a\neq \be}\la_\a\la_\be(\om_{\a\a}\otimes \om_{\be\be}+
\om_{\a\be}\otimes \om_{\be\a}).
\end{equation}

By Proposition 3.1 in \cite{XY1} we have the inequality
\begin{equation}\label{hev}
\Hess(v(\cdot,P_0))\geq v(2-v)g+\Big(\f{v-1}{pv(v^{\f{2}{p}}-1)}+\f{p+1}{pv}\Big)dv\otimes dv
\end{equation}
on $\overline \U_2$, where $g$ is the metric tensor of the $\grs{n}{m}$ and $p=min(n,m)$. Denote
\begin{equation}\label{h}
h=v^{\f{3}{2}}(2-v)^{-\f{3}{2}}
\end{equation}
on $\U_2,$ then from (4.6) in \cite{XY1}
\begin{equation}\label{Hh}
\Hess(h)=3 h\ g+\big(\f{3}{2}+\f{1}{3p}\big)h^{-1}dh\otimes dh.
\end{equation}

For a submanifold $M\to\R^{m+n}$ we have some natural functions from the functions on $\grs{n}{m}$ via the Gauss map $\g$.  We define
\begin{equation}
v\triangleq v(\cdot,P_0)\circ \g.
\end{equation}

Using the composition formula (\ref{com}), in conjunction with (\ref{hess}) and (\ref{ghij}), and the fact that  $\tau_f(\g)=0$ for translator $M$ the, we deduce the important formula as follows

\begin{equation}\label{L2v}
\mc{L}_{II} v=v|B|^2+v\sum_{\a,j}2\la_\a^2h_{\a,\a j}^2 +v\sum_{\a\neq
\be,j}\la_\a\la_\be(h_{\a,\a j}h_{\be,\be j}+h_{\a,\be
j}h_{\be,\a j})
\end{equation}

This expression can also be derived by direct computations  on the submanifold $M$, as was shown in author previous work. From the formula (3.6) of \cite{X2}
we have

\begin{equation}\label{Law}\aligned
\na_{e_i}w&=h_{\a,ij}\lan e_{\a j},\ep_1\w \cdots\w \ep_n\ran,\\
\De w=-|B|^2 w&+\sum_{\a,i,j}h_{\a,iij}\lan e_{\a i},\ep_1\w \cdots\w \ep_n\ran\\
&+ \sum_i\sum_{\a< \be,j\neq k}(h_{\a,ij} h_{\be,ik}- h_{\be,ij} h_{\a,ik})\lan e_{\a j,\be k},\ep_1\w \cdots\w \ep_n\ran.
\endaligned\end{equation}
with
\begin{equation}
e_{\a j,\be k}=e_1\w\cdots\w e_\a \w\cdots\w e_\be\w \cdots\w e_n,
\end{equation}
that is obtained by replacing $e_j$ by $e_\a$ and $e_k$ by $e_\be$ in $e_1\w \cdots\w e_n$, respectively. By (\ref{DH})
$$2^{nd}\;\text{term of} \; (\ref{Law}) =h_{\a,ij}\lan V,e_i\ran\lan e_{\a j},\ep_1\w \cdots\w \ep_n\ran=\lan V, \na w\ran.$$
Then, it follows that
\begin{equation}\label{L2w}
\mc{L}_{II} w=-|B|^2 w+ \sum_i\sum_{\a< \be,j\neq k}(h_{\a,ij} h_{\be,ik}- h_{\be,ij} h_{\a,ik})\lan e_{\a j,\be k},\ep_1\w \cdots\w \ep_n\ran.
\end{equation}

In fact, (\ref{L2v}) and (\ref{L2w}) are equivalent  each other. We summarize the following proposition.

\begin{proposition}\label{propw}
Let $M$ be a translator in $\R^{m+n}$. Then the $w-$ function satisfied (\ref{L2w}). In particular, when the normal bundle is flat (including codimension $m=1$), we have
\begin{equation}\label{L2w1}
\mc{L}_{II} w=-|B|^2 w.
\end{equation}
In general, when $m\ge 2$, if the image under the Gauss map lies in $\U$, where $w>0$, we have the formula (\ref{L2v})
\end{proposition}

Now, we  analyze   the expression (\ref{L2v}),  we group its terms according to the
different types of the indices of the coefficients of the second
fundamental form, as was done in \cite{JXY1}.  (\ref{L2v}) could be rewritten as
\begin{equation}
v^{-1}\mc{L}_{II} v=
\sum_\a\sum_{i,j>m}h_{\a,ij}^2+\sum_{j>m}I_j+\sum_{j>m,\a<\be}II_{j\a\be}
+\sum_{\a<\be<\g}III_{\a\be\g}+\sum_\a IV_\a
\end{equation} where
\begin{equation}
I_j=\sum_\a(2+2\la_\a^2)h_{\a,\a j}^2+\sum_{\a\neq
\be}\la_\a\la_\be h_{\a,\a j}h_{\be,\be j},
\end{equation}
\begin{equation}
II_{j\a\be}=2h_{\a,\be j}^2+2h_{\be,\a j}^2+2\la_\a\la_\be h_{\a,\be j}h_{\be,\a j},
\end{equation}
\begin{equation}\aligned
III_{\a\be\g}=&2h_{\a,\be\g}^2+2h_{\be,\g\a}^2+2h_{\g,\a\be}^2\\
&+2\la_\a\la_\be h_{\a,\be\g}h_{\be,\g\a}+2\la_\be\la_\g h_{\be,\g\a}h_{\g,\a\be}+2\la_\g\la_\a h_{\g,\a\be}h_{\a,\be\g}
               \endaligned
\end{equation}
and
\begin{equation}\aligned
IV_\a=&(1+2\la_\a^2)h_{\a,\a\a}^2+\sum_{\be\neq \a}\big(h_{\a,\be\be}^2+(2+2\la_\be^2)h_{\be,\be\a}^2\big)\\
      &+\sum_{\be\neq \g}\la_\be\la_\g h_{\be,\be \a}h_{\g,\g \a}+2\sum_{\be\neq \a}\la_\a\la_\be h_{\a,\be\be}h_{\be,\be\a}.
      \endaligned
\end{equation}

By  using the symmetric properties of the coefficients of the second fundamental form the the non-negative lower bounds of $I_j, II_{j\a\be},III_{\a\be\g}$
and $IV_\a$ could be estimated under the appropriate assumptions. By (\ref{Hh}) and also by using (3.13), (3.16), Lemma 3.1 and Lemma 3.2 in \cite{JXY1}, we arrive at the following
conclusion.

\begin{proposition}\label{propv}
Let $M^n$ be a translator in $\R^{n+m}$, then for arbitrary $p\in M$ and $P_0\in \grs{n}{m}$.

1. If  the image under the Gauss map of $M$ lies in $\U_2$, then from (\ref{h}) we define an auxiliary function $h$ on $\U_2$ and then  via the Gauss map $\g$ we have a function $h\circ\g$ on $M$, which is still denoted by $h$ and
\begin{equation}\label{Dh}
\mc{L}_{II}h\ge 3h|B|^2;
\end{equation}

2. If the image under the Gauss map of $M$ lies in  $\overline\U_3$ , then $\mc{L}_{II}v\ge 0$,
if the  image under the Gauss map of $M$ lies in a closed subset of $\U_3$,  namely, $v(\g(p), P_0)\le b_0<3$,  there exists a positive constant $K_0$, depending only on $b_0$, such that at $p$
\begin{equation}\label{Dv1}
\mc{L}_{II}v\geq K_0|B|^2.
\end{equation}
\end{proposition}

From (\ref{Law}) we also have
\begin{equation*}\aligned
|\na w|^2=&\sum_{\a,\be,i,j}h_{\a,ij}h_{\be,ij}A^\a A^\be\\
         =&\sum_{\a,\be}S_{\a\be}=\sum_\a\la_a|A^\a|^2\le\sum_\a\la_a=|B|^2,
         \endaligned
\end{equation*}
where $A^\a=\left<e_{\a j},\ep_1\w \cdots\w \ep_n\right>$ , $|A^\a|\le 1$ and $\la_\a$ are the eigenvalues of $S_{\a\be}$. It follows that on $\U$
\begin{equation}\label{nv}
|\na v|^2\le v^4|B|^2.
\end{equation}

\bigskip

\section{rigidity theorems}

\medskip

There is the  Sobolev inequality  \cite{MS} as follows
\begin{equation}\aligned\label{sob}
\k^{-1}\left(\int_Mg^\frac{2n}{n-2}\mathrm{d}\mu\right)^\frac{n-2}{n}\leq\int_M|\nabla
g|^2\mathrm{d}\mu+\frac{1}{2}\int_M|H|^2g^2\mathrm{d}\mu,\qquad
\forall g\in C_c^\infty(M),
\endaligned
\end{equation}
where  $\kappa>0$ is a constant.  Besides using (\ref{sob}), the
Simons type inequality of the translators \eqref{LB1} would
be used in the following result.

Let $\eta$ be a smooth function with compact support in $M$. Denote $\r=e^{\lan V, X\ran}$.
Then, we have
\begin{equation}\label{r}
\na\r=\r V^T,\quad \na\r^{\f{1}{2}}=\f{1}{2}\r^{\f{1}{2}}V^T, \quad \De\r=\r.
\end{equation}
Multiplying $\eta^2|B|^{n-2}\r$ on both sides of \eqref{LB1} and
integrating by parts yield
\begin{equation*}\aligned
0\ge&2\int_M|\na|B||^2|B|^{n-2}\e^2\rho-k\int_M|B|^{n+2}\e^2\rho -\int_M\e^2|B|^{n-2}\r\mathcal{L}_{II}|B|^2\\
=&2\int_M|\na |B||^2|B|^{n-2}\e^2\rho-k\int_M|B|^{n+2}\e^2\rho
+2\,\int_M|B|\r\na|B|\cdot\na(|B|^{n-2}\e^2)\\
=&2(n-1)\int_M|\na|B||^2|B|^{n-2}\e^2\rho-k\int_M|B|^{n+2}\e^2\rho+4\int_M(\na|B|\cdot\na\e)|B|^{n-1}\eta\rho.\\
\endaligned
\end{equation*}
By the Cauchy inequality, for any $\ep>0$, the above inequality becomes
\begin{equation}\aligned\label{3.2}
k\int_M|B|^{n+2}\e^2\rho+\frac{2}{\ep}\int_M|B|^n|\na\e|^2\r
\ge 2(n-1-\ep)\int_M|\na|B||^2|B|^{n-2}\e^2\rho.\\
\endaligned
\end{equation}
Let $f=|B|^{n/2}\r^{1/2}\e.$ Integrating by parts and using \eqref{r}, then we
have
\begin{equation}\aligned\label{2.11}
\int_M|\na f|^2=&\int_M|\na(|B|^{n/2}\e)|^2\r+\f12\int_M\na (|B|^n\e^2)\cdot\na\r+\int_M|B|^n\e^2|\na\r^{1/2}|^2\\
=&\int_M|\na(|B|^{n/2}\e)|^2\r-\f12\int_M|B|^n\e^2\De\r+\int_M|B|^n\e^2|\na\r^{1/2}|^2\\
=&\int_M|\na(|B|^{n/2}\e)|^2\r-\f 12\int_M|B|^n\e^2\r+\f{1}{4}\int_M|B|^n\e^2|V^T|^2\r.
\endaligned
\end{equation}
Combining \eqref{te}, the Sobolev inequality (\ref{sob}) and
\eqref{2.11}, we have
\begin{equation}\aligned\label{3.5}
\k^{-1}\left(\int_M|f|^{\f {2n}{n-2}}\right)^{\f{n-2}n}&\le\int_M|\na f|^2+\f12\int_M|B|^n\e^2|V^N|^2\r\\
\le&\int_M|\na(|B|^{n/2}\e)|^2\r-\f 14\int_M|B|^n\e^2\r|V^T|^2\\
=&\int_M(\f{n^2}4|\na|B||^2|B|^{n-2}\e^2+n|B|^{n-1}\e\na|B|\cdot\na\e+|B|^n|\na\e|^2)\r\\
&\hskip2in -\f 14\int_M|B|^n\e^2\r|V^T|^2.\\
\endaligned
\end{equation}
Combining the Cauchy inequality, \eqref{3.2} and \eqref{3.5},  we have
\begin{equation}\aligned\label{3.6}
&\k^{-1}\left(\int_M|f|^{\f {2n}{n-2}}\right)^{\f{n-2}n}\\
\le&\f{n^2}2\int_M|\na|B||^2|B|^{n-2}\e^2\r+2\int_M|B|^n|\na\e|^2\r-\f 14\int_M|B|^n\e^2\r|V^T|^2\\
\le&\f{n^2}{4(n-1-\ep)}\left(k\int_M|B|^{n+2}\e^2\r+\frac{2}{\ep}\int_M|B|^n|\na\e|^2\r\right)\\
&\hskip1in+2\int_M|B|^n|\na\e|^2\r-\f 14\int_M|B|^n\e^2\r|V^T|^2\\
\le&\f{kn^2}{4(n-1-\ep)}\int_M|B|^{n+2}\e^2\r+\left(\f{n^2}{2\ep(n-1-\ep)}+2\right)\int_M|B|^n|\na\e|^2\r.
\endaligned
\end{equation}
Using the H\"{o}der inequality
\begin{equation*}
\int_M|B|^{n+2}\e^2\r\le \left(\int_M|B|^{2\cdot\f n2}\right)^{\f2 n}\left(\int_M(|B|^{n}\e^2\r)^{\f n{n-2}}\right)^{\f {n-2}n}
\end{equation*}
in \eqref{3.6} yields
\begin{equation}\aligned
\k^{-1}\left(\int_M|f|^{\f {2n}{n-2}}\right)^{\f{n-2}n}
&\le\f{kn^2}{4(n-1-\ep)}\left(\int_M|B|^{2\cdot\f n2}\right)^{\f2 n}\left(\int_M(|B|^{n}\e^2\r)^{\f n{n-2}}\right)^{\f {n-2}n}\\
&\qquad+\left(\f{n^2}{2\ep(n-1-\ep)}+2\right)\int_M|B|^n|\na\e|^2\r\\
&=\f{kn^2}{4(n-1-\ep)}\left(\int_M|B|^n\right)^{\f2 n}\left(\int_M(|f|^{\f{2n}{n-2}}\right)^{\f {n-2}n}\\
&\qquad+\left(\f{n^2}{2\ep(n-1-\ep)}+2\right)\int_M|B|^n|\na\e|^2\r.
\endaligned\end{equation}
If we assume that
$$\left(\int_M|B|^n\right)^{\f{1}{n}}<\sqrt{\frac{4(n-1)}{kn^2\k}},$$
we arrive at
$$\left(\int_M|f|^{\f {2n}{n-2}}\right)^{\f{n-2}n}\le C(n,\ep)\int_M|B|^n|\na\e|^2\r,$$
where $C(n, \ep)$ is a positive constant. Choosing  an appropriate cut off function $\e$  leads the following conclusion.

\begin{theorem}
Let $M^n$ be a complete immersed translator in $\R^{n+m}.$ If
$M$ satisfies  integral conditions
$$\left(\int_M|B|^n\right)^{\f{1}{n}}<\sqrt{\frac{4(n-1)}{kn^2\k}}\qquad\text{and}\qquad \int_M|B|^ne^{\lan V, X\ran}<\infty,$$
where $\k$ is the Sobolev constant and $k=2$ when the normal bundle of $M$ in $\R^{n+m}$ is flat otherwise $k=3$, then $|B|\equiv0$ and $M$ is a linear subspace.
\end{theorem}

Let $r$ be the extrinsic  distance function from a fixed point $x_0\in M$.
Let $D_a=\{x\in M\bigcap B_a\in \R^{m+n}\}$ be closed extrinsic  ball of radius $a$ around
$x_0$ in $M$.

We now carry out  estimates for $|B|^2$ in terms of $v-$function, via auxiliary function $h$.  This local estimates implies a global
rigidity result. Such kind estimates were used in author's earlier work \cite{X3} for space-like hypersurfaces in Minkowsky space. For translator hypersurfaces
the authors \cite{BS} applied this method  and obtained a rigidity result. Now, we use this technique to higher codimension, based on the investigation of
Grassmannian manifolds.

\begin{theorem}\label{ge}
Let $M$ be a complete $n-$dimensional translator in $\R^{m+n}$ with the codimension $m\ge 2$ and the positive $w-$function.   Put $v_0=\f{2\cdot 3^{\f{2}{3}}}{1+3^{\f{2}{3}}}$. If for any constant $v_1<v_0$ the $v-$function satisfies
$$v\le v_1<v_0,$$ then $M$ is affine linear.

If for any constant $b_0<3$ the $v-$function satisfies
$$v\le b_0<3$$
and $|B|$ is constant, then $M$ is a flat Euclidean subspace.
\end{theorem}
\begin{proof}
Recall that $h=\left(\f{v}{2-v}\right)^{\f{3}{2}}$. Certainly, we have two constants
\begin{equation*}
h_0=\left(\f{v_o}{2-v_0}\right)^{\f{3}{2}}=3,\qquad h_1=\left(\f{v_1}{2-v_1}\right)^{\f{3}{2}}<3
\end{equation*}
corresponding $v_0$ and $v_1$. Choose a constant $v_2$ satisfying $v_1<v_2<v_0$. Then, we have the corresponding  constant $h_2$ satisfying $h_1<h_2=\left(\f{v_2}{2-v_2}\right)^{\f{3}{2}}<3$.
Now we begin to estimate $|B|^2,$ the squared length of the second fundamental
form of any translator $M$ in $\R^{m+n},$ in terms  the $v-$function, via $h-$function. For a fixed positive number $a,$
define $f: D_a \to \ir{}$ by
$$f=\frac {\bigl(a^2-r^2\bigr)^2|B|^2}{(h_2 - h)^2},$$
here $1\le h\le h_1<h_2$.

Since $f|_{\partial D_a}\equiv 0,\; f$ achieves an absolute maximum
in the interior of $D_a,$ say $f \le f(z),$ for some $z$ inside $D_a.$
We may also assume $|B|(z)\ne 0.$ Then from
$$\na f(z) = 0,$$
$$\Delta f(z) \le 0.$$
We obtain  the following at the point $z$:
\begin{equation}\label{df}
- \frac {2\na r^2}{a^2-r^2}+\frac {\na |B|^2}{|B|^2}+\frac {2\na h}{h_2-h}=0,
\end{equation}
\begin{equation}\label{Def}\aligned
\frac {-2|\na r^2|^2}{\bigl(a^2-r^2\bigr)^2}
-\frac {2\Delta r^2}{a^2-r^2}+\frac {\Delta |B|^2}{|B|^2}&-\frac {|\na |B|^2|^2}{|B|^4}\cr
&+\frac {2\Delta h}{h_2-h}+\frac {2|\na h|^2}{(h_2-h)^2}\le 0.
\endaligned
\end{equation}
Using the Bochner type formula (\ref{LB1}) and also (\ref{df})
\begin{equation*}\aligned
\frac {\Delta |B|^2}{|B|^2} -\frac {|\na |B|^2|^2}{|B|^4}
&\ge -\frac {2|\na h|^2}{(h_2-h)^2}
     -\frac {4|\na h||\na r^2|}{(h_2-h)(a^2-r^2)}\cr
&-\frac {2|\na r^2|^2}{(a^2-r^2)^2}- 3|B|^2-\f{1}{|B|^2}\left<V, \na|B|^2\right>.
\endaligned
\end{equation*}
Substituting the above inequality into (\ref{Def}) gives
\begin{equation}\label{ineq}\aligned
\frac {-\Delta r^2}{a^2-r^2}
&-\frac {2|\na r^2|^2}{\bigl(a^2-r^2\bigr)^2}
  -\frac {2|\na h||\na r^2|}{(h_2-h)(a^2-r^2)}\cr
&+\frac {\Delta h}{h_2-h}- \f{3}{2}|B|^2-\f{1}{2|B|^2}\left<V, \na|B|^2\right>\le 0.
\endaligned
\end{equation}

From (\ref{nv}) we have
\begin{equation}\label{nh}
|\na h|\le 3\left(\f{v_1}{2-v_1}\right)^{\f{5}{2}}|B|\triangleq C_1|B|.
\end{equation}

Considering (\ref{dr}), (\ref{Dh}) and (\ref{nh}) for estimates of $\na r^2,\; \De r^2,\De h$ and $\na h$,  (\ref{ineq}) becomes
\begin{equation}\label{ine1}\aligned
\left(\frac {3h}{h_2 - h}-\f{3}{2}\right)|B|^2
&-\left(\frac{4C_1r}{(h_2-h)(a^2-r^2)}\right)|B|\cr
&-\left(\frac {2n+2r}{a^2-r^2}+\frac {8r^2}{(a^2-r^2)^2}\right)\cr
&-\f{1}{h_2-h}\left<V, \na h\right>-\f{1}{2|B|^2}\left<V,\na|B|^2\right>\le 0.
\endaligned
\end{equation}
From (\ref{df}), we have
$$\left<\f{\na |B|^2}{2|B|^2}, V\right>+\left<\f{\na h}{h_2-h}, V\right>=\left<\f{\na r^2}{a^2-r^2}, V\right>\le \f{2r}{a^2-r^2}.$$
Therefore, the inequality (\ref{ine1}) becomes
\begin{equation}\label{ine2}\aligned
\left(\frac {3h}{ h_2 - h}-\f{3}{2}\right)|B|^2
&-\left(\frac{4C_1r}{(h_2 - h)(a^2-r^2)}\right)|B|\cr
&-\left(\frac {2n+4r}{a^2-r^2}+\frac {8r^2}{(a^2-r^2)^2}\right)\le 0.
\endaligned
\end{equation}
Note that $\f{3h}{h_2-h}-\f{3}{2}\ge\f{3}{h_2-1}-\f{3}{2}>0$. It is easily seen that if $ax^2-bx-c\le 0$ with $a,\,b,\,c$ all positive,
then
$$x^2 \le k\Bigl(\,\frac {\,b^2\,}{a^2}+ \frac {\,c\,}a\,\Bigr),$$
where $k$ is an absolute constant. In what follows $k$ may be
different in different inequalities.  Thus we obtain that, at the point $z,$
\begin{equation}
|B|^2 \le k \Biggl(\,\frac {\Bigl(\frac {4C_1r}
                        {(h_2-h)(a^2-r^2)}\Bigr)^2}
                      {\Bigl(\frac {3h}{h_2-h}
                              -\f{3}{2}\Bigr)^2}
      +\frac {(2n+4r)(a^2-r^2)+8r^2}
             {\Bigl(\frac {3h}{h_2-h}-\f{3}{2}\Bigr)
               (a^2 -r^2)^2}\,\Biggr)
\end{equation}
and
$$f(z) \le k \Biggl(\,\frac {16C_1^2r^2}
                          {\bigl(\frac {3h}{h_2-h} -\f{3}{2}\bigr)^2(h_2-h)^4}
+\frac {(2n+4r)(a^2-r^2)+8r^2}
       {\bigl(\frac {3h}{h_2-h} -\f{3}{2})(h_2-h)^2}\,\Biggr).
       $$
Since $h\ge 1$ and $h_2<3$, we have
$$\f{3h}{h_2-h}-\f{3}{2}\ge\f{3}{h_2-1}-\f{3}{2}\triangleq C_2>0, \qquad h_2-h\ge h_2-h_1.$$
Therefore,
$$f(z) \le k \Biggl(\,\frac {16C_1^2a^2}{C_2^2(h_2-h_1)^4}+\frac {(2n+4a)a^2+8a^2}{C_2(h_2-h_1)^2}\,\Biggr)$$
which implies
$$f(z)=O(a^3)\qquad\text{as}\qquad a\to\infty.$$

For any $x\in D_a,$
\begin{equation}\label{geB}
|B|^2(x)=\f{(h_2-h)^2f(x)}{(a^2-r^2)^2}\le \f{(h_2-h)^2f(z)}{(a^2-r^2)^2}\le \f{(h_2-1)^2f(z)}{(a^2-r^2)^2}.
\end{equation}
Therefore, we may fix $x$ and let $a$ tend to infinity in the above inequality. Then
we obtain $|B(x)|=0$ for all $x \in M$ and we complete the first part of the theorem.

When $v-$function satisfies
$$v\le b_0<3,$$
we could use maximum principle, Theorem \ref{max},  for the translator $M\in \R^{m+n}$.  There exists a sequence of points $\{x_k\} \in M,$ such that
$$\lim_{k\rightarrow\infty}v(x_k)=\sup f,\ \ \ \lim_{k\rightarrow\infty}|\nabla v(x_k)|=0
\ \ \ \lim_{k\rightarrow\infty}\mc{L}_{II} v(x_k)\leq0.$$
On the other hand, by Proposition \ref{propv}
$$\mc{L}_{II} v\le K_0|B|^2.$$
This means
$$\inf_M|B|^2=0,$$
which together with its constancy assumption leads our conclusion.
\end{proof}

\bibliographystyle{amsplain}

\end{document}